# ATTRACTING EDGE AND STRONGLY EDGE REINFORCED WALKS

By Vlada Limic[1] and Pierre Tarrès[2]

*CNRS, Université de Provence and Oxford University*

The goal is to show that an edge-reinforced random walk on a graph of bounded degree, with reinforcement *weight function W* taken from a general class of reciprocally summable reinforcement weight functions, traverses a random *attracting* edge at all large times.

The statement of the main theorem is very close to settling a conjecture of Sellke [Technical Report 94-26 (1994) Purdue Univ.]. An important corollary of this main result says that if $W$ is reciprocally summable and nondecreasing, the attracting edge exists on any graph of bounded degree, with probability 1. Another corollary is the main theorem of Limic [*Ann. Probab.* **31** (2003) 1615–1654], where the class of weights was restricted to reciprocally summable powers.

The proof uses martingale and other techniques developed by the authors in separate studies of edge- and vertex-reinforced walks [*Ann. Probab.* **31** (2003) 1615–1654, *Ann. Probab.* **32** (2004) 2650–2701] and of nonconvergence properties of stochastic algorithms toward unstable equilibrium points of the associated deterministic dynamics [*C. R. Acad. Sci. Sér. I Math.* **330** (2000) 125–130].

**1. Introduction.** Consider a connected graph $\mathcal{G}$ with the set of vertices $V = V(\mathcal{G})$ and the set of (unoriented) edges $E = E(\mathcal{G})$. The only assumption on the graph is that each vertex has at most $D(\mathcal{G})$ adjacent vertices (edges) for some $D(\mathcal{G}) < \infty$. So the graph $\mathcal{G}$ is either finite, or infinite with bounded degree. Special cases are infinite lattices.

Call two vertices $v, v'$ *adjacent* ($v \sim v'$ in symbols) if there exists an edge, denoted by $\{v, v'\} = \{v', v\}$, connecting them. For vertex $v$ of $\mathcal{G}$, let $A(v) \subset V$ denote the set of adjacent vertices $v' \sim v$.

Received April 2006; revised June 2006.
[1]Supported in part by an NSERC research grant and by an Alfred P. Sloan Research Fellowship.
[2]Supported in part by the Swiss National Foundation Grant 200021-1036251/1.

Let $W(k) > 0$, $k \geq 1$, be the *weight* function. The edge-reinforced random walk on $\mathcal{G}$ records a random motion of a particle along the vertices of $\mathcal{G}$ with the following properties:

  (i) if currently at vertex $v \in \mathcal{G}$, in the next step, the particle jumps to a vertex in $A(v)$;
  (ii) the probability of a jump to $v'$ is $W$-*proportional* to the number of previous traversals of the edge $\{v, v'\}$.

More precisely, let the initial *edge weights* be $X_0^e \in \mathbb{N}$ for all $e \in E$. We assume throughout the paper that $\sup_e X_0^e < \infty$.

Let $I_n$ be a $V$-valued random variable, recording the position of the particle at time $n$, $n \geq 0$. For concreteness, set $I_0 = v_0$ for some $v_0 \in V$. A *traversal* of edge $e$ occurs at time $n + 1$ if $e = \{I_n, I_{n+1}\}$. Denote by $X_n^e - X_0^e$ the total number of traversals of edge $e$ until time $n$. Let $\mathcal{F}_n$ be the filtration $\sigma\{(I_k, X_k^e, e \in E), k = 0, \ldots, n\} = \sigma\{I_k, k = 0, \ldots, n, (X_0^e, e \in E)\}$.

The *edge-reinforced random walk* on $\mathcal{G}$ with weight function $W$ is a Markov chain $(I, X) = \{(I_n, X_n^e, e \in E), n \geq 0\}$ on state space $V(\mathcal{G}) \times \mathbb{N}^E$ with the following conditional transition probabilities: on the event $\{I_n = v\}$, for $v' \in A(v)$,

$$P(I_{n+1} = v',\ X_{n+1}^e = X_n^e + \mathbb{1}_{\{e = \{v,v'\}\}},\, e \in \mathcal{G} \,|\, \mathcal{F}_n) = \frac{W(X_n^{\{v,v'\}})}{\sum_{w \in A(v)} W(X_n^{\{v,w\}})}.$$

(1)

It is easily seen that the edge-reinforced random walk is well defined for any weight function $W$, where $W(k) > 0, k \in \mathbb{N}$. Let (H0) be the following condition on $W$:

(H0) $$\sum_{k \in \mathbb{N}} \frac{1}{W(k)} < \infty.$$

Let us make a few preliminary observations. A simple calculation shows that (H0) is the necessary and sufficient condition for

$$P(\{I_n, I_{n+1}\} = \{I_0, I_1\} \text{ for all } n) > 0,$$

so that an attracting edge exists with positive probability. This implies that (H0) is necessary and sufficient to have

$$P(\text{walk is attracted to a single edge}) > 0,$$

and a variation of the above argument also implies that

$$P(\text{walk is attracted to a finite subgraph}) = 1.$$

However, it can easily be shown that if $\sum_k 1/W(k) = \infty$ and $W$ is nondecreasing, then any edge adjacent to an edge traversed by the walk infinitely often must also be traversed infinitely often.

The case $\mathcal{G} = \mathbb{Z}$ has been studied by Davis in [1], who has proven that if (H0) does not hold and $W$ is nondecreasing, then, with probability 1, every vertex in $\mathbb{Z}$ is visited by the walk infinitely often, that is, the walk is recurrent. This statement has, in fact, no proven general counterpart for other infinite graphs: even the original recurrence/transience question raised by Coppersmith and Diaconis [2] in 1986 for $W$ linear and $\mathcal{G} = \mathbb{Z}^d$, $d \geq 2$, is still open. For recent results on linearly ERRW, see [4, 6]. Sellke [7] provided examples of $W$ not nondecreasing and not satisfying the condition (H0) such that for the corresponding reinforced walk, an attractor consisting of two or more edges exists with positive (or full) probability on $\mathbb{Z}$, as well as on other graphs.

In the case where (H0) holds, the first result is due to Davis [1], who proved that on $\mathbb{Z}$, there exists almost surely some random integer $i$ such that the walk visits only $i$ and $i+1$ at all large times. Sellke [7] generalized this statement and showed that, if $\mathcal{G} = \mathbb{Z}^d$, $d \geq 1$, then under the necessary assumption (H0) (without any monotonicity requirement), there exists almost surely some random *attracting* edge, which is traversed by the walk at all large times. The same paper contains the conjecture that the above property holds for edge-reinforced random walks on general graphs of bounded degree.

The argument developed by Sellke [7] carries over to the setting where $\mathcal{G}$ is any graph of bounded degree without odd cycles, a fact used by Limic [3].

In [3], it was proven that on any graph of bounded degree, the attracting edge exists with probability 1 if $W(k) = k^\rho$ for $\rho > 1$. In this paper, we show a generalization of this result to a large class of weight functions $W$, including the class of nondecreasing weights satisfying (H0), making use of the techniques developed by the authors in [3, 8] and [9].

A *cycle* $\mathcal{C}$ in $\mathcal{G}$ of *length* $|\mathcal{C}|$ is a subgraph of $\mathcal{G}$ spanned by a $|\mathcal{C}|$-tuple of vertices $(v_1, v_2, \ldots, v_{|\mathcal{C}|})$, such that $\{v_i, v_{i+1}\} \in E$, $i = 1, \ldots, |\mathcal{C}| - 1$, $\{v_{|\mathcal{C}|}, v_1\} \in E$, and $v_i \neq v_j$ if $i \neq j$.

For each $n \in \mathbb{N}$, let

$$\alpha_n := \sum_{k \geq n} \frac{1}{W(k)^2},$$

$$\delta_n := \sum_{k=n+1}^{\infty} \left| \frac{1}{W(k)} - \frac{1}{W(k-1)} \right|.$$

Let

$$\nu(\mathcal{G}) := \sup_{\mathcal{G}' \subset \mathcal{G} \text{ odd cycle}} \sqrt{2}|\mathcal{G}'|,$$

with the convention that $\nu(\mathcal{G}) = 0$ if there are no odd cycles. Note that if there are odd cycles of arbitrarily large length in $\mathcal{G}$, then $\nu(\mathcal{G}) = \infty$.

Let (H1) be the following condition:

(H1) $$\nu(\mathcal{G}) \liminf_{n \to \infty} \frac{\delta_n}{\sqrt{\alpha_n}} < 1,$$

with the convention that $\infty \times 0 = 0 \times \infty = 0$ [i.e., if the lim inf is 0, then the condition (H1) is satisfied for any graph $\mathcal{G}$; if $\nu(\mathcal{G}) = 0$, then (H1) is satisfied for any value of the lim inf].

Let
$$\mathcal{G}_\infty = \left\{ e \in E : \sup_n X_n^e = \infty \right\}$$

be the (random) graph spanned by all edges in $\mathcal{G}$ that are traversed by the walk infinitely often. Note that

$$\{\mathcal{G}_\infty \text{ has only one edge}\} = \{\exists N < \infty \text{ s.t. } \{I_n, I_{n+1}\} = \{I_n, I_{n-1}\} \; \forall n \geq N\}$$
$$= \left\{ \exists e \in E \text{ such that } \sup_{e' \neq e} \sup_n X_n^{e'} < \infty \right\}.$$

The main result of this paper is the following theorem.

THEOREM 1. *If $W$ satisfies* (H0) *and* (H1), *then the edge-reinforced random walk on $\mathcal{G}$ traverses a random attracting edge at all large times a.s., that is,*

(2) $$P(\mathcal{G}_\infty \text{ has only one edge}) = 1.$$

Theorem 1 is proven in Section 2. It implies, in particular, that if $W$ is nondecreasing and (H0) holds, then $\mathcal{G}_\infty$ has only one edge almost surely. This statement is shown in Corollary 3 below, based on the observation that if (H1) does not hold and $W$ is nondecreasing, then $W$ belongs to a fairly large class of weights (affectionately called the *sticky weights*), given by the condition

(H2) $$\liminf_{n \to \infty} \left( \max_{0 \leq j < n} W(j) \right) \sum_{k \geq n} \frac{1}{W(k)} < \infty,$$

for which the attracting edge property is shown in Lemma 2 below. Another consequence of Theorem 1 is stated and proved at the end of this section in Corollary 4.

LEMMA 2. *(H2) implies (2).*

PROOF. Denote the finite lim inf from (H2) by $l$ and let $\mathcal{N} = \{n \geq 1 \colon (\max_{0 \leq j < n} W(j)) \sum_{k \geq n} \frac{1}{W(k)} < l + 1\}$. Then, clearly, $\mathcal{N}$ is an infinite set. Fix $m \in \mathcal{N}$ and let the ERRW $(I., X.)$ run until the time

$$J_m := \inf \left\{ k \geq 1 \colon \text{there exists an edge } e \text{ such that } X_k^e = m, \max_{e' \neq e} X_k^{e'} < m \right\},$$

at which one edge has been traversed $m$ times (where we take into account the initial number of visits $X_0^e$ to the edge) and all others have been traversed strictly fewer times. Note that the stopping time $J_m$ is a.s. finite for any $m > \sup_{e \in E} X_0^e$. Indeed, as long as all of the edges have been visited at most $m - 1$ times, the probability transitions of the random walk depend only on the values of $W$ on the set $\{1, \ldots, m - 1\}$ and are therefore bounded both below and above by positive and finite constants. Hence, for any $n < \infty$, given $\mathcal{F}_{n \wedge J_m}$ (here and later "$\wedge$" denotes the minimum operator), and on the event $\{n < J_m\}$, the probability that any particular edge $e \in E$ adjacent to the current position $I_n$ will be traversed back and forth from time $n$ until the moment its corresponding number of traversals $X_.^e$ reaches value $m$ is bounded below by a positive constant. This implies $J_m < \infty$ a.s. by the conditional Borel–Cantelli lemma.

Take $m > \sup_{e \in E} X_0^e$ and denote by $e_{J_m}$ the edge $e$ such that $X_{J_m}^e = m$. Set $d = D(\mathcal{G}) < \infty$ and note that the probability that from time $J_m$ onward the walk traverses only edge $e_{J_m}$ is bounded below by

$$
\begin{aligned}
&\prod_{k=0}^{\infty} \frac{W(m+k)}{W(m+k) + d \max_{j<m} W(j)} \\
&= \prod_{k=0}^{\infty} \left(1 - \frac{d \max_{j<m} W(j)}{W(m+k) + d \max_{j<m} W(j)}\right),
\end{aligned}
\tag{3}
$$

which is, uniformly in $m \in \mathcal{N}$, bounded away from 0 since

$$
\sum_{k=0}^{\infty} \frac{\max_{j<m} W(j)}{W(m+k) + d \max_{j<m} W(j)} < l + 1.
\tag{4}
$$

Therefore, there exists $c > 0$ such that, for all $m \in \mathcal{N}$,

$$\mathbb{E}(\mathbb{1}_{\{\text{attracting edge exists}\}} | \mathcal{F}_{J_m}) \geq c > 0.$$

Now, $(\mathcal{F}_{J_m}, m \in \mathcal{N})$ is a filtration with the natural ordering of elements of $\mathcal{N}$, and {attracting edge exists} is contained in the limiting $\sigma$-field $\lim_n \mathcal{F}_n = \lim_{m \to \infty, m \in \mathcal{N}} \mathcal{F}_{J_m}$. Here, we use the fact that $J_m$ are strictly increasing in $m$, almost surely. Therefore, the Lévy 0–1 law implies that an attracting edge must exist with probability 1. $\square$

COROLLARY 3. *Assume that $W$ is nondecreasing and that* (H0) *holds. Then $\mathcal{G}_\infty$ has one edge almost surely.*

PROOF. If $W$ is nondecreasing, then $\delta_n = 1/W(n)$. Recall that (H0) implies $\delta_n \to 0$ as $n \to \infty$. Let us prove that if $\liminf \delta_n / \sqrt{\alpha_n} > 0$, then

$$\limsup_{n \to \infty} W(n) \sum_{k \geq n} \frac{1}{W(k)} < \infty,
\tag{5}$$

implying (H2), so that an attracting edge exists almost surely by Lemma 2. This will complete the proof of the corollary since $\liminf \delta_n/\sqrt{\alpha_n} = 0$ would imply (H1).

Using the fact that $\liminf \delta_n/\sqrt{\alpha_n} > 0$, there exists $\varepsilon > 0$ such that for $n \geq n_0$,

$$\frac{1}{W(n)^2} \geq \varepsilon \sum_{k=n}^{\infty} \frac{1}{W(k)^2}.$$

This implies, for all $n \geq n_0$, that

$$\frac{1}{W(n)^2} \geq \varepsilon \sum_{k=n}^{\infty} \frac{1}{W(k)^2} \geq \varepsilon^2 \sum_{k\geq n}\sum_{j\geq k} \frac{1}{W(j)^2} = \varepsilon^2 \sum_{j\geq n} \frac{j-n+1}{W(j)^2}$$

$$\geq \varepsilon^3 \sum_{j\geq n}\sum_{k\geq j} \frac{j-n+1}{W(k)^2} \geq \frac{\varepsilon^3}{2} \sum_{k\geq n} \frac{(k-n+1)^2}{W(k)^2}.$$

Using the Cauchy–Schwarz inequality, for all $n \geq n_0$,

$$\sum_{k\geq n} \frac{1}{W(k)} = \sum_{k\geq n} \frac{k-n+1}{W(k)} \frac{1}{k-n+1}$$

$$\leq \sqrt{\sum_{k\geq n} \frac{(k-n+1)^2}{W(k)^2}} \sqrt{\sum_{k\geq n} \frac{1}{(k-n+1)^2}}$$

$$\leq \frac{\pi}{\sqrt{6}}\sqrt{\frac{2}{\varepsilon^3}}\frac{1}{W(n)},$$

which yields (5). □

REMARK. Note that no assumption on $\nu(\mathcal{G})$ is needed in the result above, nor in the next result.

Let, for all $n \geq 2$,

$$W'(n) := W(n) - W(n-1).$$

Let (H3) be the following condition:

(H3) $$\sum_{n\geq 2} \left(\frac{W'(n)}{W(n)}\right)^2 < \infty.$$

COROLLARY 4. *Assume that* (H0) *and* (H3) *hold. Then $\mathcal{G}_\infty$ has only one edge almost surely.*

PROOF. It suffices to prove that (H3) implies (H1). Suppose (H3). Then there exists $A \in \mathbb{R}_+^*$ such that for all $n \geq 2$, $W(n-1) \geq AW(n)$ [using the fact that $W(n-1)/W(n) \to 1$ as $n \to \infty$, by (H3)] and

$$\delta_n = \sum_{k \geq n} \left| \frac{1}{W(k)} - \frac{1}{W(k-1)} \right| = \sum_{k \geq n} \frac{|W'(k)|}{W(k)W(k-1)} \leq A^{-1} \sum_{k \geq n} \frac{|W'(k)|}{W(k)^2}$$

$$\leq A^{-1} \sqrt{\sum_{k \geq n} \left( \frac{W'(k)}{W(k)} \right)^2} \sqrt{\sum_{k \geq n} \frac{1}{W(k)^2}} = \sqrt{\alpha_n} \left( A^{-1} \sqrt{\sum_{k \geq n} \left( \frac{W'(k)}{W(k)} \right)^2} \right),$$

by the Cauchy–Schwarz inequality in the last line.

This last inequality yields, together with (H3), that $\limsup \delta_n / \sqrt{\alpha_n} = 0$, which implies (H1). □

REMARK. Let us give two examples of a reciprocally summable weight for which we still do not know whether or not an attracting edge exists almost surely, on graphs with at least one odd cycle. Let $W(k) := k^{1+\rho}/(2 + (-1)^k)$. Then $\delta_n \sim_{n \to \infty} 2\rho^{-1} n^{-\rho}$ and $\sqrt{\alpha_n} = O(n^{-(\rho+1/2)})$, so that (H1) is not satisfied. Assumption (H2) is not satisfied either since $\sum_{k \geq n} 1/W(k) \sim_{n \to \infty} 2\rho^{-1} n^{-\rho}$. Similarly, $W(k) = \exp\{k(2+(-1)^k)\}$, constructed from two weights satisfying hypothesis (H2) does not satisfy it anymore; nor does it satisfy (H1)—whenever $\mathcal{G}$ contains an odd cycle—since $\delta_n$ and $\sqrt{\alpha_n}$ are then of the same order asymptotically.

**2. Proof of Theorem 1.** This section is devoted to the proof of Theorem 1. The following proposition follows from results of [3] and [7].

PROPOSITION 5. *Assume that (H0) holds. Then, almost surely, $\mathcal{G}_\infty$ is either a cycle of odd length or a single edge.*

PROOF. The arguments of Section 2 in [3] apply here verbatim, but for the benefit of the reader, we provide more details. Recall that we assume throughout the paper that each vertex has at most $D(\mathcal{G})$ adjacent vertices for some $D(\mathcal{G}) < \infty$, that is, that the graph is of bounded degree. Define by $\mathcal{G}_1$ the subgraph of $\mathcal{G}$ spanned by the edges visited at least once by the walk. We know from [7], Lemma 4, that $\mathcal{G}_1$, and therefore $\mathcal{G}_\infty$, is a finite graph. First, [7], Theorem 3 (alternatively, [3], Lemma 1), shows that there is a.s. no even cycle contained in $\mathcal{G}_\infty$ and that if $\mathcal{G}_\infty$ is a tree, it a.s. only consists of two vertices and one edge connecting them. Second, [3], Lemma 2, says that there is at most one *odd* cycle contained in $\mathcal{G}_\infty$, almost surely. Third, [3], Corollary 1, says that with probability 1, $\mathcal{G}_\infty$ contains no vertex of degree $\geq 3$. Therefore, with probability 1, either $\mathcal{G}_\infty$ contains an odd cycle, in which case it is exactly equal to this cycle, or it is a single edge. □

Now, the event $\{\mathcal{G}_\infty$ is an odd cycle$\}$ is a union of at most countably many events

(6) $\qquad \{\mathcal{G}_\infty$ is an odd cycle $C\}$,

where $C$ is any fixed odd cycle in $\mathcal{G}$. Therefore, it suffices to prove that each event above happens with probability 0. Moreover, as observed in [3], if the edge-reinforced random walk on $\mathcal{G}$ stays within a finite cycle $C$ of length $\ell$ starting from some time $n_0$, at which the current edge weights on the edges of $C$ are given by $z_1, \ldots, z_\ell$, then its transition probabilities starting from time $n_0$ (and therefore its law, and asymptotic behavior) are identical to those of the edge-reinforced random walk on the cycle of length $\ell$ started from the initial configuration of edge weights $z_1, \ldots, z_\ell$. Since the event in (6) is a countable union over all finite times $n_0$ and all finite configurations $z_1, \ldots, z_\ell$, it is sufficient to prove the following proposition.

PROPOSITION 6. *Let $\mathcal{G}$ be a cycle of length $\ell$, where $\ell$ is an odd number. Assume that* (H0) *holds and that*

$$\liminf_{n \to \infty} \frac{\delta_n}{\sqrt{\alpha_n}} < \frac{1}{\sqrt{2}|\mathcal{G}|} = \frac{1}{\sqrt{2}\ell}.$$

*Then for any choice of initial condition $X_0^e \in \mathbb{N}, e \in E$, we have*

(7) $\qquad P(\mathcal{G}_\infty = \mathcal{G}) = 0.$

REMARK. For the sake of concreteness (and brevity of notation), we provide the proof of this result for the initial condition $X_0^e \equiv 1$. We remark in Section 2.3 how the proof easily extends to the general initial condition setting.

Therefore, we assume in the sequel that $\mathcal{G} := \mathbb{Z}/\ell\mathbb{Z}$, with $\ell$ odd, without loss of generality.

This section is divided into three parts: in 2.1 we introduce the processes of interest and justify them, in 2.2 we prove preliminary estimates and results, and sketch the proof of Proposition 6, which is given in 2.3.

2.1. *Preliminary notation and intuition.* For all $n \in \mathbb{N}$, let $W^*(0) = 0$ and

$$W^*(n) := \sum_{k=1}^{n} \frac{1}{W(k)}.$$

For all $n \in \mathbb{N}$ and $x \in \mathbb{Z}/\ell\mathbb{Z}$, let

$$\zeta_n(x) := \sum_{k=1}^n \left( \frac{\mathbb{1}_{\{I_{k-1}=x, I_k=x+1\}}}{W(X_{k-1}^{\{x,x+1\}})} - \frac{\mathbb{1}_{\{I_{k-1}=x+1, I_k=x\}}}{W(X_{k-1}^{\{x,x+1\}})} \right),$$

$$\varepsilon_n(x) := \sum_{k=1}^n \left( \frac{\mathbb{1}_{\{I_{k-1}=x, I_k=x+1\}}}{W(X_{k-1}^{\{x,x+1\}})} - \frac{\mathbb{1}_{\{I_{k-1}=x, I_k=x-1\}}}{W(X_{k-1}^{\{x,x-1\}})} \right),$$

$$\kappa_n(x) := \sum_{k=1}^n \left( \frac{\mathbb{1}_{\{\{I_{k-1}, I_k\}=\{x,x+1\}\}}}{W(X_{k-1}^{\{x,x+1\}})} - \frac{\mathbb{1}_{\{\{I_{k-1}, I_k\}=\{x,x-1\}\}}}{W(X_{k-1}^{\{x,x-1\}})} \right)$$

$$= W^*(X_n^{\{x,x+1\}} - 1) - W^*(X_n^{\{x,x-1\}} - 1).$$

Let us make the following observations, in order to justify the definitions of the above processes. First, note that under (H0), all of the above processes, being differences of nondecreasing bounded sequences, are bounded and have random finite limits as $n \to \infty$. Fix $x \in \mathbb{Z}/\ell\mathbb{Z}$. The process $\kappa_\cdot(x)$ is a useful way to keep track of the changes due to repeated visits of the random walk to the two edges $\{x, x-1\}$ and $\{x, x+1\}$. In particular,

(8) $$\left\{ \sup_{n \geq 1} X_n^{\{x,x-1\}} = \sup_{n \geq 1} X_n^{\{x,x+1\}} = \infty \right\} \subset \left\{ \lim_{n \to \infty} \kappa_n(x) = 0 \right\},$$

so knowing that, almost surely, $\kappa_\infty(x) = \lim_{n \to \infty} \kappa_n(x) \neq 0$ for at least one $x \in \mathbb{Z}/\ell\mathbb{Z}$ would be sufficient to conclude (7). The proof of Proposition 6 relies on this observation.

The process $\varepsilon_\cdot(x)$ is analytically the nicest of the three since it is a martingale.

LEMMA 7. *For each $x \in \mathbb{Z}/\ell\mathbb{Z}$, the process $(\varepsilon_n(x))_{n \in \mathbb{N}}$ is a martingale.*

PROOF. Given $x \in \mathbb{Z}/\ell\mathbb{Z}$, $(\varepsilon_n(x))_{n \in \mathbb{N}}$ is a martingale since for all $n \in \mathbb{N}$, $\varepsilon_{n+1}(x) = \varepsilon_n(x)$ if $I_n \neq x$ and if $I_n = x$, then

$$\mathbb{E}(\varepsilon_{n+1}(x) - \varepsilon_n(x) | \mathcal{F}_n)$$

$$= \frac{W(X_n^{\{x,x+1\}})}{W(X_n^{\{x,x+1\}}) + W(X_n^{\{x,x-1\}})} \frac{1}{W(X_n^{\{x,x+1\}})}$$

$$- \frac{W(X_n^{\{x,x-1\}})}{W(X_n^{\{x,x-1\}}) + W(X_n^{\{x,x+1\}})} \frac{1}{W(X_n^{\{x,x-1\}})} = 0. \qquad \square$$

Note that $\varepsilon_\cdot(x)$ only captures half of the traversals of edges $\{x, x+1\}$ and $\{x, x-1\}$, namely those originating from the central vertex $x$.

Process $\zeta_\cdot(x)$ is a measure of difference in the directional visits to edge $\{x, x+1\}$. Clearly,

$$\sum_{k=1}^{n} \frac{\mathbb{1}_{\{\{I_{k-1}, I_k\} = \{x, x+1\}\}}}{W(X_{k-1}^{\{x,x+1\}})} = 2 \sum_{k=1}^{n} \frac{\mathbb{1}_{\{I_{k-1}=x, I_k=x+1\}}}{W(X_{k-1}^{\{x,x+1\}})} - \zeta_n(x)$$

and, similarly,

$$\sum_{k=1}^{n} \frac{\mathbb{1}_{\{\{I_{k-1}, I_k\} = \{x, x-1\}\}}}{W(X_{k-1}^{\{x,x-1\}})} = 2 \sum_{k=1}^{n} \frac{\mathbb{1}_{\{I_{k-1}=x, I_k=x-1\}}}{W(X_{k-1}^{\{x,x-1\}})} + \zeta_n(x-1).$$

A useful relation follows:

$$(9) \qquad \kappa_n(x) = 2\varepsilon_n(x) - \zeta_n(x) - \zeta_n(x-1).$$

Moreover, note that for all $n \in \mathbb{N}$,

$$(10) \qquad \sum_{x \in \mathbb{Z}/\ell\mathbb{Z}} (\zeta_n(x) - \varepsilon_n(x)) = 0$$

since

$$\zeta_n(x) - \varepsilon_n(x) = \sum_{k=1}^{n} \left( \frac{\mathbb{1}_{\{I_{k-1}=x, I_k=x-1\}}}{W(X_{k-1}^{\{x,x-1\}})} - \frac{\mathbb{1}_{\{I_{k-1}=x+1, I_k=x\}}}{W(X_{k-1}^{\{x,x+1\}})} \right),$$

which implies that

$$\sum_{x \in \mathbb{Z}/\ell\mathbb{Z}} (\zeta_n(x) - \varepsilon_n(x))$$

$$= \sum_{x \in \mathbb{Z}/l\mathbb{Z}} \sum_{k=1}^{n} \frac{\mathbb{1}_{\{I_{k-1}=x, I_k=x-1\}}}{W(X_{k-1}^{\{x,x-1\}})} - \sum_{x \in \mathbb{Z}/l\mathbb{Z}} \sum_{k=1}^{n} \frac{\mathbb{1}_{\{I_{k-1}=x+1, I_k=x\}}}{W(X_{k-1}^{\{x+1,x\}})}$$

$$= \sum_{k=1}^{n} \sum_{x \in \mathbb{Z}/l\mathbb{Z}} \frac{\mathbb{1}_{\{I_{k-1}=x, I_k=x-1\}}}{W(X_{k-1}^{\{x,x-1\}})} - \sum_{k=1}^{n} \sum_{x \in \mathbb{Z}/l\mathbb{Z}} \frac{\mathbb{1}_{\{I_{k-1}=x, I_k=x-1\}}}{W(X_{k-1}^{\{x,x-1\}})} = 0.$$

Recall that $\delta_n = \sum_{k=n+1}^{\infty} |\frac{1}{W(k)} - \frac{1}{W(k-1)}|$. For all $k, n \in \mathbb{N}$ such that $n \geq k$, let

$$\delta_{k,n}(x) := \delta_{X_k^{\{x,x+1\}}} - \delta_{X_n^{\{x,x+1\}}}, \qquad \text{if } k \leq n < \infty,$$

$$\delta_{k,\infty}(x) := \delta_{X_k^{\{x,x+1\}}}, \qquad \text{if } k \leq \infty,$$

$$\Delta_{k,n} := \sum_{x \in \mathbb{Z}/\ell\mathbb{Z}} \delta_{k,n}(x), \qquad \text{if } k \leq n \leq \infty.$$

Note that

$$(11) \qquad \frac{1}{W(X_k^{\{x,x+1\}})} \leq \sum_{j=X_k^{\{x,x+1\}}}^{\infty} \left| \frac{1}{W(j)} - \frac{1}{W(j+1)} \right| = \delta_{k,\infty}.$$

Fix $m \in \mathbb{N}$. For all $n \in \mathbb{N} \cup \{\infty\}$ and $x \in \mathbb{Z}/\ell\mathbb{Z}$, let $X_n^x$ be the number of times the *vertex* $x$ has been visited during time interval $[m, n]$:

$$X_n^x := \sum_{k=m}^{n} \mathbb{1}_{\{I_k = x\}}.$$

REMARK. Recall that for $e \in \mathcal{E}$, $X_n^e$ is the number of times plus $X_0^e$ that *edge* $e$ has been visited up to and including time $n$. The new notation will not cause confusion since edges will be always denoted either by letters $e, f$ or sets $\{\cdot, \cdot\}$.

For each $n \in \mathbb{N} \cup \{\infty\}$ and $x \in \mathbb{Z}/\ell\mathbb{Z}$, let $t_n(x)$ (that also depends on $m$ fixed above) be the time of $n$th visit to $x$ during interval $[m, \infty]$:

$$t_n(x) := \inf\{k \geq m : X_k^x = n\} = \inf\{k > t_{n-1}(x) : I_k = x\}.$$

Note that $t_n(x)$ may take the value $\infty$ (if $X_k^x < n, \forall k$) and then $t_j(x) = \infty$ for all $j \geq n$. However, if $t_n(x) < \infty$, then $t_{n+1}(x) > t_n(x)$, almost surely.

In the proof, we will focus on one particular vertex of the cycle, adjacent to the least-visited edge at some particular time. We will suppose it is vertex 0 for simplicity and let

$$\kappa_n := \kappa_n(0), \qquad t_n := t_n(0).$$

Note that

(12)
$$\begin{aligned} P(\mathcal{G}_\infty \neq \mathcal{G} | \mathcal{F}_m) \\ &= \mathbb{E}(P(\mathcal{G}_\infty \neq \mathcal{G} | \mathcal{F}_{t_n}) | \mathcal{F}_m) \\ &= \mathbb{E}(P(\mathcal{G}_\infty \neq \mathcal{G} | \mathcal{F}_{t_n}) \mathbb{1}_{\{t_n = \infty\}} + P(\mathcal{G}_\infty \neq \mathcal{G} | \mathcal{F}_{t_n}) \mathbb{1}_{\{t_n < \infty\}} | \mathcal{F}_m) \\ &\geq P(t_n = \infty | \mathcal{F}_m) + \mathbb{E}(\mathbb{E}(\mathbb{1}_{\{\mathcal{G}_\infty \neq \mathcal{G}\}} | \mathcal{F}_{t_n}) \mathbb{1}_{\{t_n < \infty\}} | \mathcal{F}_m). \end{aligned}$$

We will provide a uniform lower bound

(13) $$\mathbb{E}(\mathbb{1}_{\{\mathcal{G}_\infty \neq \mathcal{G}\}} | \mathcal{F}_{t_n}) \mathbb{1}_{\{t_n < \infty\}} \geq c \mathbb{1}_{\{t_n < \infty\}}$$

for some $c > 0$, for all $n \in \mathbb{N}$ sufficiently large [cf. discussion preceding (26)]. This will imply that

$$P(\mathcal{G}_\infty \neq \mathcal{G} | \mathcal{F}_m) \geq P(t_n = \infty | \mathcal{F}_m) + (1 - P(t_n = \infty | \mathcal{F}_m))c \geq c$$

for all $m$ sufficiently large, and the Lévy 0–1 law will imply (7).

2.2. *Preliminary results.* The proof of Proposition 6 is based on a study of the behavior of $\kappa_n(x)$, $x \in \mathbb{Z}/\ell\mathbb{Z}$. We prove in Section 2.3 that there exists a.s. $x \in \mathbb{Z}/\ell\mathbb{Z}$ such that $\kappa_n(x)$ does not converge to 0, which enables one to conclude that, almost surely, either $\{x, x+1\}$ or $\{x, x-1\}$ is visited finitely often [otherwise, $\kappa_\infty(x) = 0$].

Given $x \in \mathbb{Z}/\ell\mathbb{Z}$, the process $(\kappa_n(x))_{n \geq 0}$, contrary to $(\varepsilon_n(x))_{n \geq 0}$, is not a martingale. Our first aim is therefore to estimate its mean behavior. Note that if $\ell$ were even, we would be able to answer the question without such an estimate, by the construction of a martingale $(R_n(x))_{n \geq 0}$ combining the processes $(\kappa_n(x))_{n \geq 0}$, $x \in \mathbb{Z}/\ell\mathbb{Z}$:

$$R_n(x) := \sum_{x \in \mathbb{Z}/\ell\mathbb{Z}, x \text{ even}} \kappa_n(x) = \sum_{x \in \mathbb{Z}/\ell\mathbb{Z}} (-1)^x \varepsilon_n(x).$$

Then an upper bound of the variance of the increments would enable us to prove that $R_n(x)$ a.s. does not converge to 0, which subsequently implies that there is at least one $x \in \mathbb{Z}/\ell\mathbb{Z}$ such that $\kappa_n(x)$ does not converge to 0, as required. Sellke [7] obtains the corresponding result (for $\ell$ even) using a construction due to Rubin.

The behavior of $(\kappa_n(x))_{n \in \mathbb{N}}$ is described by equation (9):

$$\kappa_n(x) = 2\varepsilon_n(x) - \zeta_n(x) - \zeta_n(x-1),$$

where $\zeta_n(y)$, $y \in \mathbb{Z}/l\mathbb{Z}$, defined in Section 2.1, is the difference between the weighted numbers of visits from $y$ to $y+1$ and those from $y+1$ to $y$. Hence, the study of $\kappa_n(x)$ requires [through $\zeta_n(x)$ and $\zeta_n(x-1)$] some information on the probabilities of cycles

$$\vec{q}_n := P(\{I_{t_n+1} = 1\} \cap \{I_{t_{n+1}-1} = -1 \text{ or } t_{n+1} = \infty\}|\mathcal{F}_{t_n}),$$
$$\overleftarrow{q}_n := P(\{I_{t_n+1} = -1\} \cap \{I_{t_{n+1}-1} = 1 \text{ or } t_{n+1} = \infty\}|\mathcal{F}_{t_n}).$$

The quantity $\vec{q}_n$ defined here is the probability of a *cycle from the right* (from 0 to 1 and then returning 0 by $-1$), whereas $\overleftarrow{q}_n$ is the probability of a *cycle from the left* (from 0 to $-1$ and then return to 0 by 1).

A natural method to compare $\vec{q}_n$ and $\overleftarrow{q}_n$ would be to write them explicitly for each possible path and then to find a coupling of the corresponding paths with reversed paths, as done in [3] for $W(n) = n^\rho$, $\rho > 1$. But this method is difficult to apply since the estimates depend on the regularity of $W$, as well as on the current numbers of visits to the edges.

We relate the quantities $\vec{q}_n$ and $\overleftarrow{q}_n$ using the two following observations:

(i) for any $x \in \mathbb{Z}/\ell\mathbb{Z}$, $\mathbb{E}(\zeta_{t_{n+1}}(x) - \zeta_{t_n}(x)|\mathcal{F}_{t_n})$ provides a good estimate of $(\vec{q}_n - \overleftarrow{q}_n)/W(X_{t_n}^{\{x,x+1\}})$, provided $\mathbb{E}(\delta_{t_n,t_{n+1}}(x)|\mathcal{F}_{t_n})$ is small [see Lemma 8(i) and (ii)];

(ii) $\sum_{x \in \mathbb{Z}/\ell\mathbb{Z}} \zeta_n(x)$ is a martingale [Lemma 7 and equation (10)].

Therefore, since all of these estimates of $\mathbb{E}(\zeta_{t_{n+1}}(x) - \zeta_{t_n}(x)|\mathcal{F}_{t_n})$ have the sign of $\vec{q}_n - \overleftarrow{q}_n$ and sum to zero, $\vec{q}_n - \overleftarrow{q}_n$ is negligible and $(\kappa_{t_n})_{n\in\mathbb{N}}$ is close to a martingale with respect to filtration $(\mathcal{F}_{t_n})_{n\in\mathbb{N}}$ [see Lemma 8(iv)]. Here, we also need the fact that under assumption (H1), $\mathbb{E}(\delta_{t_n,t_{n+1}}(x)|\mathcal{F}_{t_n})$ can be neglected for any $x \in \mathbb{Z}/\ell\mathbb{Z}$.

The link between $\mathbb{E}(\zeta_{t_{n+1}}(x) - \zeta_{t_n}(x)|\mathcal{F}_{t_n})$ and $(\vec{q}_n - \overleftarrow{q}_n)/W(X_{t_n}^{\{x,x+1\}})$ described above is a consequence of the fact that the evolution $\zeta_{t_{n+1}}(x) - \zeta_{t_n}(x)$ is only significant over the excursions $(t_n, t_{n+1})$ away from 0 which are cycles, where it increases (resp., decreases) by $1/W(X_{t_n}^{\{x,x+1\}})$ if the cycle is from the right (resp., from the left), while during the excursions which are not cycles (this happens whenever $I_{t_n+1} = I_{t_{n+1}-1}$), for each $x$, the traversals of an edge $\{x, x+1\}$ contribute as many times positively as negatively to the evolution of $\zeta.(x)$.

The property that $(\kappa_{t_n})_{n\in\mathbb{N}}$ is close to a martingale enables one to control the evolution of $\kappa_{t_n}^2$ (Lemma 10) and to prove (in Section 2.3) that $\kappa_{t_n}$ does not converge to 0 with lower bounded probability if $\{0,1\}$ is at some point traversed fewer times than the other edges $\{x, x+1\}$, $x \neq 0$.

For all $x \in \mathbb{Z}/\ell\mathbb{Z}$, let us define the $(\mathcal{F}_{t_n})_{n\geq 2}$-adapted processes $(u_n(x))_{n\geq 2}$, $(\Lambda_n)_{n\geq 2}$, $(\lambda_n(x))_{n\geq 2}$ and $(v_n)_{n\geq 2}$ by

$$u_n(x) := \zeta_{t_n}(x) - \zeta_{t_{n-1}}(x)$$

(14)
$$- \frac{\mathbb{1}_{\{I_{t_{n-1}+1}=1\}\cap\{I_{t_n-1}=-1 \text{ or } t_n=\infty\}} - \mathbb{1}_{\{I_{t_{n-1}+1}=-1\}\cap\{I_{t_n-1}=1 \text{ or } t_n=\infty\}}}{W(X_{t_{n-1}}^{\{x,x+1\}})},$$

$$\Lambda_n := \sum_{x\in\mathbb{Z}/\ell\mathbb{Z}} \frac{1}{W(X_{t_n}^{\{x,x+1\}})},$$

(15)
$$\lambda_n(x) := \frac{1/W(X_{t_n}^{\{x,x+1\}})}{\Lambda_n} \in (0,1)$$

and

$$v_n := \sum_{x\in\mathbb{Z}/\ell\mathbb{Z}} (\lambda_{n-1}(0) + \lambda_{n-1}(-1) - \mathbb{1}_{x\in\{0,-1\}})u_n(x).$$

The processes $(u_n(x))_{n\geq 2}$ and $(v_n)_{n\geq 2}$ play an important role, as made explicit in the following lemma: $\mathbb{E}(v_{n+1}|\mathcal{F}_{t_n})$ is the drift increment of $\kappa_{t_{n+1}} - \kappa_{t_n}$ [Lemma 8(iii)] and $u_n(x)$ and $v_{n+1}$ [which is a weighted sum of $u_{n+1}(x)$, $x \in \mathbb{Z}/\ell\mathbb{Z}$] are small, by Lemma 8(ii), $|u_n(x)| \leq \delta_{t_{n-1},t_n}(x)$. Recall that we concentrate here on the behavior of the process $\kappa_n \equiv \kappa_n(0)$.

LEMMA 8. *For all $x \in \mathbb{Z}/\ell\mathbb{Z}$ and $n \in \mathbb{N}$, a.s. on $\{t_n < \infty\}$:*

(i) $\mathbb{E}(\zeta_{t_{n+1}}(x) - \zeta_{t_n}(x) - u_{n+1}(x)|\mathcal{F}_{t_n}) = \frac{\overrightarrow{q_n} - \overleftarrow{q_n}}{W(X_{t_n}^{\{x,x+1\}})}$,

(ii) $|u_{n+1}(x)| \le \delta_{t_n,t_{n+1}}(x), |v_{n+1}| \le \Delta_{t_n,t_{n+1}}$,

(iii) $\mathbb{E}(\kappa_{t_{n+1}} - \kappa_{t_n} - v_{n+1}|\mathcal{F}_{t_n}) = 0$,

(iv) $|\mathbb{E}(\kappa_{t_{n+1}} - \kappa_{t_n}|\mathcal{F}_{t_n})| \le \mathbb{E}(\Delta_{t_n,t_{n+1}}|\mathcal{F}_{t_n})$.

PROOF. Property (i) follows directly from definition (14).

Let us prove property (ii): Assume that $t_n < \infty$ and note that by symmetry, it suffices to consider the case $I_{t_n+1} = 1$. We then have

$$\zeta_{t_{n+1}}(x) - \zeta_{t_n}(x) = \sum_{k=X_{t_n}^{\{x,x+1\}}}^{X_{t_{n+1}}^{\{x,x+1\}} - 1} \frac{(-1)^{k - X_{t_n}^{\{x,x+1\}}}}{W(k)} \tag{16}$$

since during the time interval $(t_n, t_{n+1})$, there are $X_{t_{n+1}}^{\{x,x+1\}} - X_{t_n}^{\{x,x+1\}}$ (possibly infinitely many) traversals of the edge $\{x, x+1\}$ in alternating directions and since the first traversal (if there is one) happens in the direction of the directed edge $(x, x+1)$.

Assume first that, in addition, $t_{n+1} < \infty$. Now, either $I_{t_{n+1}-1} = I_{t_n+1} = 1$ or $I_{t_{n+1}-1} = -1$. In the former case, there is an even number of summands in (16) with alternating signs for each $x$, and, clearly,

$$|u_{n+1}(x)| = |\zeta_{t_{n+1}}(x) - \zeta_{t_n}(x)| \le \delta_{t_n,t_{n+1}}(x).$$

In the latter case, (16) consists of an odd number of terms and, similarly,

$$|u_{n+1}(x)| = |\zeta_{t_{n+1}}(x) - \zeta_{t_n}(x) - 1/W(X_{t_n}^{\{x,x+1\}})| \le \delta_{t_n,t_{n+1}}(x). \tag{17}$$

Next, assume that $t_{n+1} = \infty$. Then there exists $y \ne 0$ such that $\{y, y+1\}$ becomes the attracting edge during the "uncompleted excursion" $[t_n, \infty)$. The reasoning is very similar to the one above. Namely, if $x < y$, then the sum in (16) consists of an odd number of alternating terms and, again, an estimate (17) applies. If $x > y$, then (16) has an even number of terms and we use estimate (11) to derive (17). If $x = y$, then (16) is an infinite alternating sum, so the reasoning from the $t_{n+1} < \infty$ case applies.

To bound $v_{n+1}$, use the fact that $\lambda_n(0) + \lambda_n(-1) \in [0,1]$ and conclude that

$$|v_{n+1}| \le \sum_{x \in \mathbb{Z}/l\mathbb{Z}} |u_{n+1}(x)| \le \sum_{x \in \mathbb{Z}/l\mathbb{Z}} \delta_{t_n,t_{n+1}}(x) = \Delta_{t_n,t_{n+1}}.$$

Let us now prove (iii): for all $x \in \mathbb{Z}/\ell\mathbb{Z}$, $(\varepsilon_n(x))_{n \in \mathbb{N}}$ are martingales and equation (10) holds, therefore $(\sum_{x \in \mathbb{Z}/\ell\mathbb{Z}} \zeta_n(x))_{n \in \mathbb{N}}$ is a martingale. This implies, summing (i) over $x \in \mathbb{Z}/\ell\mathbb{Z}$, that

$$(\overleftarrow{q_n} - \overrightarrow{q_n}) \sum_{x \in \mathbb{Z}/\ell\mathbb{Z}} \frac{1}{W(X_{t_n}^{\{x,x+1\}})} = \sum_{x \in \mathbb{Z}/\ell\mathbb{Z}} \mathbb{E}(u_{n+1}(x)|\mathcal{F}_{t_n})$$

or, alternatively,

$$\overleftarrow{q_n} - \overrightarrow{q_n} = \frac{\mathbb{E}(\sum_{x \in \mathbb{Z}/\ell\mathbb{Z}} u_{n+1}(x)|\mathcal{F}_{t_n})}{\Lambda_n}. \tag{18}$$

Hence, property (i) for $x \in \mathbb{Z}/\ell\mathbb{Z}$ implies that on $\{t_n < \infty\}$, we have

$$\mathbb{E}(\zeta_{t_{n+1}}(x) - \zeta_{t_n}(x)|\mathcal{F}_{t_n}) = \mathbb{E}(u_{n+1}(x)|\mathcal{F}_{t_n}) - \lambda_n(x) \sum_{y \in \mathbb{Z}/\ell\mathbb{Z}} \mathbb{E}(u_{n+1}(y)|\mathcal{F}_{t_n})$$

$$= \sum_y (\mathbb{1}_{\{y=x\}} - \lambda_n(x))\mathbb{E}(u_{n+1}(y)|\mathcal{F}_{t_n}).$$

Accordingly, using (9) and Lemma 7, we conclude that on $\{t_n < \infty\}$,

$$\mathbb{E}(\kappa_{t_{n+1}} - \kappa_{t_n}|\mathcal{F}_{t_n}) = - \sum_{x \in \{0,-1\}} \mathbb{E}(\zeta_{t_{n+1}}(x) - \zeta_{t_n}(x)|\mathcal{F}_{t_n})$$

$$= \sum_{x \in \mathbb{Z}/\ell\mathbb{Z}} (\lambda_n(0) + \lambda_n(-1) - \mathbb{1}_{x \in \{0,-1\}})\mathbb{E}(u_{n+1}(x)|\mathcal{F}_{t_n})$$

$$= \mathbb{E}[v_{n+1}|\mathcal{F}_{t_n}].$$

Property (iv) follows from (ii) and (iii). □

Even more precise estimates of the drift of $\kappa.$ will be needed and the following technical lemma provides the necessary calculations.

LEMMA 9. (i) *On* $\{t_n < \infty\} \cap \{I_{t_n+1} = 1\}$, *we have*

$$\mathbb{E}(\kappa_{t_{n+1}} - \kappa_{t_n+1} - v_{n+1}|\mathcal{F}_{t_n+1}) = \lambda_n(0)(1 - \lambda_n(0) - \lambda_n(-1))\Lambda_n.$$

(ii) *On* $\{t_n < \infty\} \cap \{I_{t_n+1} = -1\}$, *we have*

$$\mathbb{E}(\kappa_{t_{n+1}} - \kappa_{t_n+1} - v_{n+1}|\mathcal{F}_{t_n+1}) = -\lambda_n(-1)(1 - \lambda_n(0) - \lambda_n(-1))\Lambda_n.$$

REMARK. Note that Lemma 8(iii) is a consequence of (i)–(ii) by nested conditioning on $\mathcal{F}_{t_n+1}$, although we preferred to give its proof independently for the benefit of the reader.

PROOF OF LEMMA 9. Let us assume that $t_n < \infty$ and $I_{t_n+1} = 1$, and prove (i).

Let us define

$$\overrightarrow{r_n} := P(I_{t_{n+1}-1} = -1 \text{ or } t_{n+1} = \infty|\mathcal{F}_{t_n+1}),$$

$$\overleftarrow{r_n} := P(I_{t_{n+1}-1} = 1 \text{ or } t_{n+1} = \infty|\mathcal{F}_{t_n+1}).$$

Then by (14), for all $x \in \mathbb{Z}/\ell\mathbb{Z} \setminus \{0\}$, using the fact that $\zeta_{t_n}(x) = \zeta_{t_n+1}(x)$ and $I_{t_n+1} = 1$, we have

$$\mathbb{E}(\zeta_{t_{n+1}}(x) - \zeta_{t_n+1}(x)|\mathcal{F}_{t_n+1}) = \mathbb{E}(\zeta_{t_{n+1}}(x) - \zeta_{t_n}(x)|\mathcal{F}_{t_n+1})$$

$$= \mathbb{E}(u_{n+1}(x)|\mathcal{F}_{t_n+1}) + \frac{\vec{r_n}}{W(X_{t_n}^{\{x,x+1\}})}$$

and, similarly using $\zeta_{t_n+1}(0) = \zeta_{t_n}(0) + 1/W(X_{t_n}^{\{0,1\}})$,

$$\mathbb{E}(\zeta_{t_{n+1}}(0) - \zeta_{t_n+1}(0)|\mathcal{F}_{t_n+1}) = \mathbb{E}(u_{n+1}(0)|\mathcal{F}_{t_n+1}) + \frac{\vec{r_n}-1}{W(X_{t_n}^{\{0,1\}})}.$$

Again using the fact that $(\sum_{x \in \mathbb{Z}/\ell\mathbb{Z}} \zeta_n(x))_{n \in \mathbb{N}}$ is a martingale, we obtain

$$\sum_{x \in \mathbb{Z}/\ell\mathbb{Z}} \mathbb{E}(\zeta_{t_{n+1}}(x) - \zeta_{t_n+1}(x)|\mathcal{F}_{t_n+1})$$

$$= 0 = \vec{r_n}\Lambda_n + \left( \sum_{x \in \mathbb{Z}/\ell\mathbb{Z}} \mathbb{E}(u_{n+1}(x)|\mathcal{F}_{t_n+1}) - \frac{1}{W(X_{t_n}^{\{0,1\}})} \right).$$

Therefore,

$$\vec{r_n} = \left( \lambda_n(0) - \frac{1}{\Lambda_n} \sum_{x \in \mathbb{Z}/\ell\mathbb{Z}} \mathbb{E}(u_{n+1}(x)|\mathcal{F}_{t_n+1}) \right),$$

which, by (9), implies that

$$\mathbb{E}(\kappa_{t_{n+1}} - \kappa_{t_n+1}|\mathcal{F}_{t_n+1})$$

$$= - \sum_{x \in \{0,-1\}} \mathbb{E}(\zeta_{t_{n+1}}(x) - \zeta_{t_n+1}(x)|\mathcal{F}_{t_n+1})$$

$$= \frac{1}{W(X_{t_n}^{\{0,1\}})} - \vec{r_n}\left( \frac{1}{W(X_{t_n}^{\{0,1\}})} + \frac{1}{W(X_{t_n}^{\{0,-1\}})} \right)$$

$$\quad - \sum_{x \in \{0,-1\}} \mathbb{E}(u_{n+1}(x)|\mathcal{F}_{t_n+1})$$

$$= \lambda_n(0)\Lambda_n - \left( \lambda_n(0) - \frac{1}{\Lambda_n} \sum_{x \in \mathbb{Z}/\ell\mathbb{Z}} \mathbb{E}(u_{n+1}(x)|\mathcal{F}_{t_n+1}) \right)(\lambda_n(0) + \lambda_n(-1))\Lambda_n$$

$$\quad - \sum_{x \in \{0,-1\}} \mathbb{E}(u_{n+1}(x)|\mathcal{F}_{t_n+1})$$

$$= \sum_{x \in \mathbb{Z}/\ell\mathbb{Z}} (\lambda_n(0) + \lambda_n(-1) - \mathbb{1}_{x \in \{0,-1\}})\mathbb{E}(u_{n+1}(x)|\mathcal{F}_{t_n+1})$$

$$+ \lambda_n(0)\Lambda_n(1 - \lambda_n(0) - \lambda_n(-1))$$
$$= \mathbb{E}[v_{n+1}|\mathcal{F}_{t_n+1}] + \lambda_n(0)\Lambda_n(1 - \lambda_n(0) - \lambda_n(-1)).$$

One can similarly show (ii). $\square$

LEMMA 10. *For all $n \in \mathbb{N}$, a.s. on $\{t_n < \infty\}$,*

$$\mathbb{E}(\kappa_{t_{n+1}}^2 - \kappa_{t_n}^2|\mathcal{F}_{t_n}) \geq \mathbb{E}\left(\sum_{k=X_{t_n}^{\{0,1\}}}^{X_{t_{n+1}}^{\{0,1\}}-1} \frac{1}{W(k)^2} + \sum_{k=X_{t_n}^{\{0,-1\}}}^{X_{t_{n+1}}^{\{0,-1\}}-1} \frac{1}{W(k)^2}\bigg|\mathcal{F}_{t_n}\right)$$
$$- 2\mathbb{E}(|\kappa_{t_n+1}|\Delta_{t_n,t_{n+1}}|\mathcal{F}_{t_n}).$$

PROOF. We split
$$\kappa_{t_{n+1}}^2 - \kappa_{t_n}^2 = (\kappa_{t_{n+1}}^2 - \kappa_{t_n+1}^2) + (\kappa_{t_n+1}^2 - \kappa_{t_n}^2)$$
and compute the conditional expectation of each summand separately. First, $\kappa_{t_n+1} - \kappa_{t_n} = \varepsilon_{t_n+1} - \varepsilon_{t_n}$ implies that

(19) $$\mathbb{E}(\kappa_{t_n+1} - \kappa_{t_n}|\mathcal{F}_{t_n}) = 0$$

and, therefore,

(20)
$$\mathbb{E}(\kappa_{t_n+1}^2 - \kappa_{t_n}^2|\mathcal{F}_{t_n}) = \mathbb{E}((\kappa_{t_n+1} - \kappa_{t_n})^2|\mathcal{F}_{t_n})$$
$$= \mathbb{E}\left(\frac{\mathbb{1}_{\{I_{t_n+1}=1\}}}{W(X_{t_n}^{\{0,1\}})^2} + \frac{\mathbb{1}_{\{I_{t_n+1}=-1\}}}{W(X_{t_n}^{\{0,-1\}})^2}\bigg|\mathcal{F}_{t_n}\right).$$

Next, we wish to estimate $\mathbb{E}(\kappa_{t_{n+1}}^2 - \kappa_{t_n+1}^2|\mathcal{F}_{t_n})$ from below. Note first that

(21)
$$\kappa_{t_{n+1}}^2 - \kappa_{t_n+1}^2$$
$$= (\kappa_{t_{n+1}} - \kappa_{t_n+1})^2 + 2\kappa_{t_n+1}(\kappa_{t_{n+1}} - \kappa_{t_n+1})$$
$$= (\kappa_{t_{n+1}} - \kappa_{t_n+1})^2 + 2\kappa_{t_n+1}(\kappa_{t_{n+1}} - \kappa_{t_n+1} - v_{n+1}) + 2\kappa_{t_n+1}v_{n+1}$$
$$= (\kappa_{t_{n+1}} - \kappa_{t_n+1})^2 + 2\kappa_{t_n}(\kappa_{t_{n+1}} - \kappa_{t_n+1} - v_{n+1})$$
$$+ 2(\kappa_{t_n+1} - \kappa_{t_n})(\kappa_{t_{n+1}} - \kappa_{t_n+1} - v_{n+1}) + 2\kappa_{t_n+1}v_{n+1}.$$

Now,

(22)
$$\mathbb{E}((\kappa_{t_{n+1}} - \kappa_{t_n+1})^2|\mathcal{F}_{t_n})$$
$$= \mathbb{E}\left(\frac{\mathbb{1}_{\{I_{t_{n+1}-1}=1, t_{n+1}<\infty\}}}{W(X_{t_{n+1}-1}^{\{0,1\}})^2} + \frac{\mathbb{1}_{\{I_{t_{n+1}-1}=-1, t_{n+1}<\infty\}}}{W(X_{t_{n+1}-1}^{\{0,-1\}})^2}\bigg|\mathcal{F}_{t_n}\right).$$

Lemma 8(iii) and identity (19) imply that

(23) $$\mathbb{E}(2\kappa_{t_n}(\kappa_{t_{n+1}} - \kappa_{t_n+1} - v_{n+1})|\mathcal{F}_{t_n}) = 0.$$

Lemma 9(i)–(ii) implies that $\mathbb{E}(\kappa_{t_{n+1}} - \kappa_{t_n+1} - v_{n+1}|\mathcal{F}_{t_n+1})$ is positive when $I_{t_n+1} = 1$ and negative when $I_{t_n+1} = -1$, hence that it has the same sign as $\kappa_{t_n+1} - \kappa_{t_n}$. Therefore, by nested conditioning with respect to $\mathcal{F}_{t_n+1}$, we obtain that

(24) $$\mathbb{E}(2(\kappa_{t_n+1} - \kappa_{t_n})(\kappa_{t_{n+1}} - \kappa_{t_n+1} - v_{n+1})|\mathcal{F}_{t_n}) \geq 0.$$

Adding together inequalities (20), (22), (23) and (24) completes the proof, using the second inequality in Lemma 8(ii). □

2.3. *Proof of Proposition* 6. Assume that the conditions of Proposition 6 hold. Fix $a < 1$, for which the subset of $\mathbb{N}$ defined by

(25) $$\Theta \equiv \Theta_a := \left\{ n \in \mathbb{N} : \delta_n \leq \frac{a\sqrt{\alpha_n}}{\sqrt{2\ell}} \right\} \qquad \text{is infinite.}$$

Let

$$\mathcal{A} := \left\{ \min_{x \in \mathbb{Z}/\ell\mathbb{Z}} X_\infty^{\{x,x+1\}} < \infty \right\} = \{\mathcal{G}_\infty \neq \mathcal{G}\}.$$

As remarked earlier, it suffices to prove that there exists a constant $C \in \mathbb{R}_+^*$ such that, for all $m \in \mathbb{N}$ sufficiently large, $P(\mathcal{A}|\mathcal{F}_m) \geq C$.

Let $m \in \mathbb{N}$. Let, for all $n \in \mathbb{N}$,

$$o_n \in \underset{e \in E(\mathbb{Z}/\ell\mathbb{Z})}{\arg\min} \{X_n^e\}$$

so that $o_n$ is an edge (if there is more than one edge minimizing $X_n^\cdot$, choose $o_n$ arbitrarily from the set of minima) corresponding to the smallest $X_n^\cdot$. Note that if $X_0^e \equiv 1$, then $o_n$ corresponds to the least-visited edge at time $n$.

Define stopping time $k_1 \equiv k_1(m) := \inf\{n > m : X_n^{o_n} \in \Theta \setminus \{X_m^{o_m}\}\}$. Without loss of generality, assume that $\{k_1 < \infty\}$ happens, since

(26) $$P(\mathcal{A}|\mathcal{F}_m) = \mathbb{E}(\mathbb{1}_\mathcal{A}|\mathcal{F}_m) = \mathbb{E}(\mathbb{1}_\mathcal{A}\mathbb{1}_{\{k_1=\infty\}}|\mathcal{F}_m) + \mathbb{E}(\mathbb{1}_\mathcal{A}\mathbb{1}_{\{k_1<\infty\}}|\mathcal{F}_m)$$
$$= P(\{k_1 = \infty\}|\mathcal{F}_m) + \mathbb{E}(P(\mathcal{A}|\mathcal{F}_{k_1})\mathbb{1}_{\{k_1<\infty\}}|\mathcal{F}_m),$$

where we use the fact that $\Theta$ is infinite and therefore $\{k_1 = \infty\} \subset \{\mathcal{G} \neq \mathcal{G}_\infty\} = \mathcal{A}$.

Without loss of generality, assume that $o_{k_1} = \{0, 1\}$ and $I_{k_1} = 0$, and note that we have $t_{k_0} = k_1$, for some (random) positive integer $k_0$.

Therefore, it suffices to find a positive lower bound on $P(\mathcal{A}|\mathcal{F}_{k_1}) = P(\mathcal{A}|\mathcal{F}_{t_{k_0}})$, on the event $\{k_1 < \infty\}$.

Let
$$\varepsilon := 1 \wedge \frac{a^{-1}-1}{1+\sqrt{2}}.$$

For all $n \in \mathbb{N}$, define
$$e(n) := X_{t_n}^{\{0,1\}} \wedge X_{t_n}^{\{0,-1\}}$$
to be the number of traversals of the *weaker* (i.e., visited fewer times) edge at 0. Note, in particular, that $e(k_0) = \min\{X_{t_{k_0}}^{\{0,1\}}, X_{t_{k_0}}^{\{0,-1\}}\} = X_{t_{k_0}}^{\{0,1\}} = X_{k_1}^{\{0,1\}} \in \Theta$, and that $n \mapsto e(n)$ is nondecreasing, which implies, in particular, that $n \mapsto \alpha_{e(n)}$ and $n \mapsto \delta_{e(n)}$ are nonincreasing. We next state a few similar and easy facts for future reference. For all $n \geq k_0$,

(27) $\quad \delta_{e(n)} \leq \delta_{e(k_0)}, \qquad 2\Delta_{t_n,t_{n+1}} \leq 2\Delta_{t_{k_0},\infty} \leq 2\ell\delta_{e(k_0)} \leq a\sqrt{2\alpha_{e(k_0)}}.$

Define the $(\mathcal{F}_{t_n})_{n \geq k_0}$ stopping time
$$S := \inf\{n \geq k_0 : |\kappa_{t_n}| \geq \varepsilon\sqrt{\alpha_{e(n)}} + \Delta_{t_n,\infty}\}$$
$$\wedge \inf\{n > k_0 : |\kappa_{t_{n-1}+1}| \geq \varepsilon\sqrt{\alpha_{e(k_0)}} + \Delta_{t_{k_0},\infty}\}.$$

The next three lemmas make use of the techniques developed in Lemmas 1 and 2 in [8], and in [5].

REMARK. To explain the technique informally, consider the martingale $M_{\cdot} = \kappa_{t_{\cdot}} - \mathrm{drift}(\kappa_{t_{\cdot}})$. Then the first lemma says that $|M_{\cdot}|$ infinitely often becomes larger than a fixed small $\varepsilon$ proportion of the total standard deviation $\sqrt{\alpha_{e(\cdot)}} = \mathrm{SD}(M_{\infty} - M_{\cdot})$ of the (infinitely many) remaining increments of $M$; the next two lemmas say that in each situation above, there is positive (bounded away from 0) probability that $|M_{\cdot}|$ remains strictly above value 0, either by not exiting $(\varepsilon\sqrt{\alpha_{e(\cdot)}}/2, 4\sqrt{\alpha_{e(\cdot)}})$ or by exiting it through the larger boundary point and not coming back to 0 due to Doob's $L^2$ inequality. However, if $\kappa_n \to 0$ as $n \to \infty$ then $M_n \to 0$ as $n \to \infty$.

REMARK. We prove Proposition 6 for the initial condition $X_0^e \equiv 1, e \in E$, but the same proof carries through for a general initial condition, after a few minor modifications: the first line in the definition of the $\kappa_{\cdot}(x)$ processes changes and the second line remains the same,
$$\kappa_n(x) = W^*(X_n^{\{x,x+1\}} - 1) - W^*(X_n^{\{x,x-1\}} - 1),$$
so that the important relations (9)–(10) are replaced by an equality up to a constant. This does not modify the results in Section 2.2, which only use estimates of differences $\kappa_{n+1} - \kappa_n$ and $\zeta_{n+1} - \zeta_n$. The goal is still to prove that $\kappa_n$ does not converge to 0 a.s.

LEMMA 11. *For $a$ and $\varepsilon$ defined above, $P(\{S < \infty\} \cup \mathcal{A}|\mathcal{F}_{t_{k_0}}) \geq a\varepsilon/19$.*

PROOF. Assume that $S > k_0$ and let, for $n > k_0$,
$$z_n := \kappa_{t_n}^2 + \sqrt{2\alpha_{e(k_0)}}(|\kappa_{t_{n-1}+1}| - \varepsilon\sqrt{\alpha_{e(k_0)}} - \Delta_{t_{k_0},\infty})^+,$$
where $x^+ = \max(x,0)$. Then on $\{n < S\}$,
$$z_{n+1} - z_n = \kappa_{t_{n+1}}^2 - \kappa_{t_n}^2 + \sqrt{2\alpha_{e(k_0)}}(|\kappa_{t_n+1}| - \varepsilon\sqrt{\alpha_{e(k_0)}} - \Delta_{t_{k_0},\infty})^+$$
and using (27), it is easy to check that
$$\begin{aligned}
(28) \quad 2\Delta_{t_n,t_{n+1}}|\kappa_{t_n+1}| &\leq 2\Delta_{t_n,t_{n+1}}(\varepsilon\sqrt{\alpha_{e(k_0)}} + \Delta_{t_{k_0},\infty}) \\
&\quad + \sqrt{2\alpha_{e(k_0)}}(|\kappa_{t_n+1}| - \varepsilon\sqrt{\alpha_{e(k_0)}} - \Delta_{t_{k_0},\infty})^+.
\end{aligned}$$

Due to (28) and Lemma 10, on $\{k_0 \leq n < S\}$,
$$\begin{aligned}
\mathbb{E}(z_{n+1} - z_n|\mathcal{F}_{t_n}) &\geq \mathbb{E}\left(\sum_{k=X_{t_n}^{\{0,1\}}}^{X_{t_{n+1}}^{\{0,1\}}-1} \frac{1}{W(k)^2} + \sum_{k=X_{t_n}^{\{0,-1\}}}^{X_{t_{n+1}}^{\{0,-1\}}-1} \frac{1}{W(k)^2} \bigg| \mathcal{F}_{t_n}\right) \\
&\quad - 2(\varepsilon\sqrt{\alpha_{e(k_0)}} + \Delta_{t_{k_0},\infty})\mathbb{E}(\Delta_{t_n,t_{n+1}}|\mathcal{F}_{t_n}).
\end{aligned}$$

A careful reader will note that the definition of stopping time $S$ is designed precisely to give the inequality above, that is, to eliminate the auxiliary term $-2|\kappa_{t_n+1}|\Delta_{t_n,t_{n+1}}$ in the drift estimate for $\kappa_{t_n}^2$ in Lemma 10.

Therefore, by uniform integrability,
$$\begin{aligned}
&\mathbb{E}(z_S - z_{k_0}|\mathcal{F}_{t_{k_0}}) \\
&\geq \mathbb{E}\left(\sum_{k=X_{t_{k_0}}^{\{0,1\}}}^{X_{t_{S-1}}^{\{0,1\}}-1} \frac{1}{W(k)^2} + \sum_{k=X_{t_{k_0}}^{\{0,-1\}}}^{X_{t_{S-1}}^{\{0,-1\}}-1} \frac{1}{W(k)^2} \bigg| \mathcal{F}_{t_{k_0}}\right) \\
(29) \quad &\quad - 2(\varepsilon\sqrt{\alpha_{e(k_0)}} + \Delta_{t_{k_0},\infty})\mathbb{E}(\Delta_{t_{k_0},t_\infty}|\mathcal{F}_{t_{k_0}}) \\
&\geq P(\{S = \infty\} \cap \mathcal{A}^c|\mathcal{F}_{t_{k_0}})\alpha_{e(k_0)} - 2(\varepsilon\sqrt{\alpha_{e(k_0)}} + \ell\delta_{e(k_0)})\mathbb{E}(\Delta_{t_{k_0},t_\infty}|\mathcal{F}_{t_{k_0}}) \\
&\geq [P(\{S = \infty\} \cap \mathcal{A}^c|\mathcal{F}_{t_{k_0}}) - (1 - a\varepsilon)]\alpha_{e(k_0)},
\end{aligned}$$
using, in the third inequality, the fact that $e(k_0) \in \Theta$ [cf. (25)], so
$$\begin{aligned}
2(\varepsilon\sqrt{\alpha_{e(k_0)}} + \ell\delta_{e(k_0)})\Delta_{t_{k_0},t_\infty} &\leq 2(\varepsilon\sqrt{\alpha_{e(k_0)}} + \ell\delta_{e(k_0)})\ell\delta_{e(k_0)} \\
&\leq \sqrt{2}a\varepsilon\alpha_{e(k_0)} + a^2\alpha_{e(k_0)} \\
&\leq a(1 + \sqrt{2}\varepsilon)\alpha_{e(k_0)} \leq (1 - a\varepsilon)\alpha_{e(k_0)},
\end{aligned}$$

where the last inequality follows from the definition of $\varepsilon$.

On the other hand, note that $(\kappa_{t_{n+1}} - \kappa_{t_n})^2 \leq 4\alpha_{e(n)}$ for all $n$ and if $n \in [k_0, S)$ then,

$$\kappa_{t_{n+1}}^2 = (\kappa_{t_n} + \kappa_{t_{n+1}} - \kappa_{t_n})^2 \leq 2(\kappa_{t_n}^2 + 4\alpha_{e(n)})$$

$$\leq 2\left(\left(\frac{a}{\sqrt{2}} + \varepsilon\right)^2 \alpha_{e(k_0)} + 4\alpha_{e(k_0)}\right) \leq 16\alpha_{e(k_0)}$$

and $\kappa_{t_n+1} \leq \kappa_{t_n} + \sqrt{\alpha_{e(n)}} \leq \kappa_{t_n} + \sqrt{\alpha_{e(k_0)}}$. Hence, on $\{S < \infty\}$, we have

$$z_S \leq (16 + \sqrt{2})\alpha_{e(k_0)} \leq 18\alpha_{e(k_0)},$$

on $\{S = \infty\}$, we have

$$z_S \leq \left(\varepsilon + \frac{a}{\sqrt{2}}\right)^2 \alpha_{e(k_0)} \leq 18\alpha_{e(k_0)}$$

and on $\{S = \infty\} \cap \mathcal{A}^c$, we have $z_S = 0$, so

(30) $$\mathbb{E}(z_S - z_{k_0}|\mathcal{F}_{t_{k_0}}) \leq 18\alpha_{e(k_0)}P(\{S < \infty\} \cup \mathcal{A}|\mathcal{F}_{t_{k_0}}),$$

since $z_{k_0} \geq 0$.

If we combine (29) with (30) and the two preceding inequalities, we obtain, for $p := P(\{S < \infty\} \cup \mathcal{A}|\mathcal{F}_{t_{k_0}})$,

$$18p \geq 1 - p - (1 - a\varepsilon) = a\varepsilon - p,$$

which implies the lemma. $\square$

Now, assume that $S < \infty$ and, for instance, $\kappa_{t_{S-1}+1} \geq \varepsilon\sqrt{\alpha_{e(k_0)}} + \Delta_{t_{k_0},\infty}$ with $S > k_0$, and define the $(\mathcal{F}_{t_n})$ stopping time

$$U := \inf\left\{n \geq S : \kappa_{t_n} \notin \left(\frac{\varepsilon}{2}\sqrt{\alpha_{e(k_0)}} + \Delta_{t_n,\infty}, 4\sqrt{\alpha_{e(k_0)}} + \Delta_{t_n,\infty}\right)\right\}.$$

REMARK. The two remaining cases where $S$ happens due to $|\kappa_{t_S}| > \varepsilon\sqrt{\alpha_{e(S)}} + \Delta_{t_S,\infty}$, can be dealt with in a very similar way. Assuming that $\kappa_{t_S} > 0$, one would redefine

$$U := \inf\left\{n > S : \kappa_{t_n} \notin \left(\frac{\varepsilon}{2}\sqrt{\alpha_{e(S)}}/2 + \Delta_{t_n,\infty}, 4\sqrt{\alpha_{e(S)}} + \Delta_{t_n,\infty}\right)\right\}$$

and, in the statement (and proof) of the next lemma, $e(k_0)$ would have to be replaced by $e(S)$ and $t_{S-1}+1$ by $t_S$ everywhere, and the estimate Lemma 9(i) would no longer be needed.

Note that if $U = \infty$, then it is not the case that both $\{0,1\}$ and $\{0,-1\}$ are visited infinitely often, so $\mathcal{A}$ happens.

LEMMA 12. *If $S > k_0$ and $\kappa_{t_{S-1}+1} \geq \varepsilon\sqrt{\alpha_{e(k_0)}} + \Delta_{t_{k_0},\infty}$, then*

$$P(\{U = \infty\} \cup \{U < \infty, \kappa_{t_U} \geq 4\sqrt{\alpha_{e(k_0)}} + \Delta_{t_U,\infty}\}|\mathcal{F}_{t_{S-1}+1}) \geq \varepsilon/16.$$

PROOF. Using Lemma 8(iii) with $n \in [S, U)$, Lemma 9(i) with $n = S - 1$, and nested expectation, we get

$$\text{(31)} \qquad \mathbb{E}\left(\kappa_{t_U} - \kappa_{t_{S-1}+1} - \sum_{k=S-1}^{U-1} v_{k+1} \Big| \mathcal{F}_{t_{S-1}+1}\right) \geq 0.$$

Now, let us assume that $S > k_0$ and $\kappa_{t_{S-1}+1} \geq \varepsilon\sqrt{\alpha_{e(k_0)}} + \Delta_{t_{k_0},\infty}$, and let us define

$$W := \kappa_{t_U} - \Delta_{t_U,\infty} - (\kappa_{t_{S-1}+1} - \Delta_{t_{k_0},\infty}) \leq \kappa_{t_U} - \Delta_{t_U,\infty} - \varepsilon\sqrt{\alpha_{e(k_0)}}.$$

Inequality (31) implies, using Lemma 8(ii), that

$$\text{(32)} \qquad \mathbb{E}(W) = \mathbb{E}(\kappa_{t_U} - \kappa_{t_{S-1}+1} + \Delta_{t_{k_0},t_U}) \geq 0.$$

Now, note that on $\{U < \infty, \kappa_{t_U} \leq \varepsilon\sqrt{\alpha_{e(k_0)}}/2 + \Delta_{t_U,\infty}\}$, we have, by definition,

$$\text{(33)} \qquad W \leq \kappa_{t_U} - \Delta_{t_U,\infty} - \varepsilon\sqrt{\alpha_{e(k_0)}} \leq -\varepsilon\sqrt{\alpha_{e(k_0)}}/2.$$

On the other hand, on the complement $\{U = \infty\} \cup \{U < \infty, \kappa_{t_U} \geq 4\sqrt{\alpha_{e(k_0)}} + \Delta_{t_U,\infty}\}$ (and in fact on the whole probability space),

$$\text{(34)} \qquad W \leq \kappa_{t_U} - \Delta_{t_U,\infty} \leq 7\sqrt{\alpha_{e(k_0)}}.$$

Indeed, if $U = \infty$,

$$\kappa_{t_U} - \Delta_{t_U,\infty} = \lim_{n \to \infty}(\kappa_{t_n} - \Delta_{t_n,\infty}) \leq 4\sqrt{\alpha_{e(k_0)}},$$

by definition of $U$, and if $U < \infty$, then, using (27),

$$\kappa_{t_U} - \Delta_{t_U,\infty} = \kappa_{t_{U-1}} - \Delta_{t_{U-1},\infty} + (\kappa_{t_U} - \kappa_{t_{U-1}})$$
$$+ (\Delta_{t_{U-1},\infty} - \Delta_{t_U,\infty})$$
$$\leq \kappa_{t_{U-1}} - \Delta_{t_{U-1},\infty} + 3\sqrt{\alpha_{e(U-1)}} \leq 7\sqrt{\alpha_{e(k_0)}}.$$

In summary, equations (32), (33) and (34) imply, letting

$$p := P(\{U = \infty\} \cup \{U < \infty, \kappa_{t_U} \geq 4\sqrt{\alpha_{e(k_0)}} + \Delta_{t_U,\infty}\}|\mathcal{F}_{t_{S-1}+1}),$$

that

$$\text{(35)} \qquad 0 \leq \mathbb{E}(W) \leq -(1-p)\varepsilon\sqrt{\alpha_{e(k_0)}}/2 + 7\sqrt{\alpha_{e(k_0)}}p,$$

which completes the proof. □

LEMMA 13. *If $U < \infty$ and $\kappa_{t_U} \geq 4\sqrt{\alpha_{e(U)}} + \Delta_{t_U,\infty}$, then*

$$P(\mathcal{A}|\mathcal{F}_{t_U}) \geq P\left(\liminf_{n\to\infty} \kappa_{t_n} > 0 \Big| \mathcal{F}_{t_U}\right) \geq 7/16.$$

PROOF. Define the $(\mathcal{F}_{t_n})$-adapted process

$$\xi_n := \kappa_{t_n} - \sum_{k \leq n} v_k.$$

By Lemma 8(iii), $(\xi_n)_{n\geq 1}$ is a bounded martingale (note that this is not a Doob–Meyer decomposition), and hence converges a.s.; this can also be seen as a consequence of the convergences of $\kappa_{t_n}$ and $\sum v_k$, being differences of nondecreasing bounded sequences.

Since $\sum_{k>U} |v_k| \leq \Delta_{t_U,\infty}$ by Lemma 8(ii), we have

(36) $\quad P(\kappa_{t_\infty} = \liminf \kappa_{t_n} > 0 | \mathcal{F}_{t_U}) \geq P(|\xi_\infty - \xi_U| < 4\sqrt{\alpha_{e(U)}} | \mathcal{F}_{t_U}).$

Now, due to a martingale property of $\xi$,

(37) $$\mathbb{E}((\xi_\infty - \xi_U)^2 | \mathcal{F}_{t_U}) = \mathbb{E}\left(\sum_{k=U}^{\infty} (\xi_{k+1} - \xi_k)^2 \Big| \mathcal{F}_{t_U}\right).$$

To estimate the right-hand side of (37), note that for all $k \geq U$, using

$$\xi_{k+1} - \xi_k = (\kappa_{t_{k+1}} - \kappa_{t_k+1}) + (\kappa_{t_k+1} - \kappa_{t_k}) + v_{k+1},$$

we obtain that

$$(\xi_{k+1} - \xi_k)^2 \leq 3[(\kappa_{t_{k+1}} - \kappa_{t_k+1})^2 + (\kappa_{t_k+1} - \kappa_{t_k})^2 + v_{k+1}^2]$$
$$\leq 3[(\kappa_{t_{k+1}} - \kappa_{t_k+1})^2 + (\kappa_{t_k+1} - \kappa_{t_k})^2 + \Delta_{t_k,t_{k+1}}^2].$$

This implies, using identities (20) and (22), and nested expectation, that

(38)
$$\mathbb{E}\left(\sum_{k=U}^{\infty} (\xi_{k+1} - \xi_k)^2 \Big| \mathcal{F}_{t_U}\right)$$
$$= \mathbb{E}\left(\sum_{k=U}^{\infty} \mathbb{E}((\xi_{k+1} - \xi_k)^2 | \mathcal{F}_{t_k}) \Big| \mathcal{F}_{t_U}\right)$$
$$\leq 3\mathbb{E}\left(2 \sum_{j=X_{t_U}^{\{0,1\}} \wedge X_{t_U}^{\{0,-1\}}}^{\infty} \frac{1}{W(k)^2} + \Delta_{t_U,\infty}^2 \Big| \mathcal{F}_{t_U}\right) \leq 9\alpha_{e(U)}.$$

Using (37) and (38), together with the Markov inequality, we obtain

$$P(|\xi_\infty - \xi_U| \geq 4\sqrt{\alpha_{e(U)}} | \mathcal{F}_{t_U}) \leq \frac{\mathbb{E}((\xi_\infty - \xi_U)^2 | \mathcal{F}_{t_U})}{16\alpha_{e(U)}} \leq \frac{9}{16},$$

which gives the conclusion, by (36). □

The three lemmas above enable us to complete the proof. Indeed, let us define $\tilde{t_S}$ by $\tilde{t_S} := t_{S-1} + 1$ if $S > k_0$ and $\kappa_{t_{S-1}+1} \geq \varepsilon\sqrt{\alpha_{e(k_0)}} + \Delta_{t_{k_0},\infty}$, and $\tilde{t_S} = t_S$ otherwise. Then $\{S < \infty\}$ is $\mathcal{F}_{\tilde{t_S}}$-measurable. By Lemmas 12 and 13 (and the remark preceding Lemma 12), the same argument as in equality (26) yields, using $\{U = \infty\} \subset \mathcal{A}$, that if $S < \infty$, then

$$\begin{aligned}P(\mathcal{A}|\mathcal{F}_{\tilde{t_S}}) \\ &= P(U = \infty|\mathcal{F}_{\tilde{t_S}}) + \mathbb{E}(P(\mathcal{A}|\mathcal{F}_{t_U})\mathbb{1}_{\{U<\infty\}}|\mathcal{F}_{\tilde{t_S}}) \\ &\geq P(U = \infty|\mathcal{F}_{\tilde{t_S}}) + \mathbb{E}(P(\mathcal{A}|\mathcal{F}_{t_U})\mathbb{1}_{\{U<\infty, \kappa_{t_U} \geq 4\sqrt{\alpha_{e(k_0)}} + \Delta_{t_U,\infty}\}}|\mathcal{F}_{\tilde{t_S}}) \\ &\geq P(\{U = \infty\} \cup \{U < \infty, \kappa_{t_U} \geq 4\sqrt{\alpha_{e(k_0)}} + \Delta_{t_U,\infty}\}|\mathcal{F}_{\tilde{t_S}}) \geq 7\varepsilon/256.\end{aligned}$$

Now, this inequality yields, together with Lemma 11, that

$$\begin{aligned}P(\mathcal{A}|\mathcal{F}_{t_{k_0}}) &\geq \mathbb{E}(\mathbb{1}_\mathcal{A}\mathbb{1}_{\{S=\infty\}}|\mathcal{F}_{t_{k_0}}) + \mathbb{E}(\mathbb{1}_{\{S<\infty\}}P(\mathcal{A}|\mathcal{F}_{\tilde{t_S}})|\mathcal{F}_{t_{k_0}}) \\ &\geq P(\mathcal{A} \cup \{S < \infty\}|\mathcal{F}_{t_{k_0}}) \times 7\varepsilon/256 \geq 7a\varepsilon^2/4864.\end{aligned}$$

Due to inequality (26) and the Lévy 0–1 law, we conclude that Proposition 6 holds.

**Acknowledgments.** P. T. wishes to thank the Department of Mathematics at the University of British Columbia for their hospitality during a month-long visit in the summer of 2004. V. L. wishes to thank the Laboratoire de Statistique et Probabilités at Université Paul Sabatier for their hospitality during a visit in the summer of 2003.

We are grateful to the anonymous referee for very useful comments and suggestions.

Centre de Mathématiques et Informatique  
Université de Provence, LATP UMR 6632  
Technopôle Château-Gombert  
39, rue F. Joliot Curie  
13453 Marseille Cedex 13  
France  
E-mail: vlada@cmi.univ-mrs.fr

Mathematical Institute  
University of Oxford  
24–29 St Giles'  
Oxford OX1 3LB  
United Kingdom  
E-mail: tarres@maths.ox.ac.uk